\documentclass[10pt,a4paper,twoside]{article}

\usepackage[english]{babel}
\usepackage{array,multirow}
\usepackage{a4wide} 
\usepackage{subcaption} 
\usepackage{epsfig} 
\usepackage{amsmath,amsthm,amsfonts,amssymb} 
\usepackage{float} 
\usepackage{url} 
\usepackage{array}
\usepackage{algorithm}
\usepackage{algorithmic}
\usepackage{graphicx}
\usepackage{epstopdf}
\usepackage{cite}

\newtheorem{Theorem}{Theorem}[section]
\newtheorem{definition}[Theorem]{Definition}
\newtheorem{proposition}[Theorem]{Proposition}

\newtheorem{Lemma}[Theorem]{Lemma}

\newcommand\blfootnote[1]{%
	\begingroup
	\renewcommand\thefootnote{}\footnote{#1}%
	\addtocounter{footnote}{-1}%
	\endgroup
}

\begin{document}

\title{Boundary Controllability of Riemann--Liouville\\ 
Fractional Semilinear Equations\blfootnote{This 
is a preprint version of the paper published open access in 
'Commun. Nonlinear Sci. Numer. Simul.' [https://doi.org/10.1016/j.cnsns.2023.107814].
Submitted 26/Jul/2022; Revised 08/Dec/2022 and 16/Oct/2023;
Accepted for publication 30/Dec/2023; Available online 03/Jan/2024.}}

\author{Asmae Tajani$^{1,2}$\\
\url{https://orcid.org/0000-0002-3631-5376}\\
{\tt tajaniasmae@ua.pt}
\and Fatima-Zahrae El Alaoui$^{2}$\\
\url{https://orcid.org/0000-0001-8912-4031}\\
{\tt f.elalaoui@umi.ac.ma} 
\and Delfim F. M. Torres$^{1,}$\thanks{Corresponding author:
telephone: +351 234 370 668; fax: +351 234 370 066; e-mail: delfim@ua.pt}\\ 
\url{https://orcid.org/0000-0001-8641-2505}\\
{\tt delfim@ua.pt}}

\date{$^{1}$\mbox{Center for Research and Development in Mathematics and Applications (CIDMA),} 
Department of Mathematics, University of Aveiro, 3810-193 Aveiro, Portugal\\[0.3cm]
$^{2}$TSI Team, Department of Mathematics, Faculty of Sciences,\\
Moulay Ismail University, 11201 Meknes, Morocco}


\maketitle


\begin{abstract}
We study the boundary regional controllability of a class 
of Riemann--Liouville fractional semilinear sub-diffusion 
systems with boundary Neumann conditions. The result is obtained 
by using semi-group theory, the fractional Hilbert uniqueness method, 
and Schauder's fixed point theorem. Conditions on the order 
of the derivative, internal region, and on the nonlinear part 
are obtained. Furthermore, we present appropriate sufficient 
conditions for the considered fractional system to be regionally 
controllable and, therefore, boundary regionally controllable. 
An example of a population density system with diffusion is given 
to illustrate the obtained theoretical results. Numerical simulations 
show that the proposed method provides satisfying results 
regarding two cases of the control operator.

\medskip


\end{abstract}
   

\section{Introduction}

Mathematical control theory may be seen as an interdisciplinary branch 
of engineering and mathematics that deals with the behavior of dynamic 
systems through many perceptions, one of them being controllability.
In short, controllability means that one can govern the system to any 
desired state of the dynamical system, in finite time, using a suitable 
control function into the system. The need to control a system has been 
visualized in many phenomena, such as the design of electric circuits, 
manufacturing processes, communication systems, and biological networks 
\cite{ref3,ref1,ref2,bio2}. Therefore, the controllability question 
has been widely studied by many authors, for linear and nonlinear systems,  
in finite and Hilbert state spaces, etc. 
(see, e.g., \cite{cc2,delay,edo,nolin} and references therein). 

In many phenomena, we are interested to govern the system into 
a desired state in a subregion of the whole domain and, in some cases, 
we are also interested to impose some computational requirements. 
For this reason, an increasing of interest related to regional controllability 
can be observed in the literature \cite{reg1,reg3,MR4342539,reg2}.  

Regional and approximate controllability of partial differential equations (PDEs) 
has been the object of intensive researches during the last decade
\cite{MR3244467,MR3895586,Torres2018,Torres2020}.
In particular, fractional PDEs become more and more interesting 
for modeling anomalous phenomena in complex systems theory. 
Indeed, memory effects, that appear in some complex processes,  
can be better described by using time-fractional systems   
than using integer-order equations \cite{MR4200529,MR3999702,MR4232864}.
This is in particular true for population models, where the fractional 
operators make the increasing behavior of a population slower. This has 
been well documented, for instance, by Mirzazadeh et al., who presented 
a nonlinear time fractional biological population model and established 
an efficient method to find its exact solution \cite{pop}. 

Due to the large number of applications of fractional-order systems, 
several papers on controllability of fractional control systems 
have appeared in the literature \cite{ref6,ref5,apcont,contsl}. 
For the regional controllability problem, Ge et al. studied 
the regional controllability of a linear diffusion system 
with a Caputo derivative of order $\alpha\in(0,1]$ 
in two cases: when the control operator is bounded, 
by using the Hilbert Uniqueness Method (HUM) approach; 
and when one has an unbounded operator, 
by using a compactness condition \cite{art1}. 
Tajani et al. established sufficient conditions 
for regional controllability of Riemann--Liouville fractional 
semilinear sub-diffusion systems using two approaches: 
an analytical approach, based on the contraction mapping theorem; 
and using the HUM approach, where the regional approximate controllability 
of the associated linear system is assumed to hold \cite{mee,ma}.\\

Our goal is to establish the regional boundary 
controllability, that is, to provide a suitable control
that is able to bring the state of the system to a desired one 
in the boundary of the whole domain. That is the case, 
for example, of a biological reactor in which the concentration 
regulation of a substratum at the bottom of the reactor is expected. 
A fractional system of Riemann--Liouville type can effectively capture 
the complex behavior of the biological reactor, where the fractional 
derivative represents the asymmetric heat transfer and the hysteresis 
effect of temperature variation. Tusset et al. 
studied the controllability of a cooling fluid flow that passes 
through the fermented jacket, in order to maintain the ideal temperature 
in the biological reactor in two cases: integer-order system and 
Riemann--Liouville type fractional system, in addition to some successful 
numerical results in both cases \cite{bio}. To do the regional boundary 
controllability, our approach consists to study the regional controllability 
in an internal subregion that contains a part of the boundary, 
based on a generalization of the Gronwall--Bellman inequality,
and then to bring the state to the boundary by projection.

Motivated by the above arguments, and to see the deviation from 
the classical behavior, which is studied numerically for some real phenomena, 
for example the cholera outbreak studied in \cite{bal}, 
when Baleanu et al. demonstrate that the fractional model of order approximated 1 
is more consistent with the real data when compared to the classical model. 
Our new techniques allow us to establish the regional boundary 
controllability of a wider class of fractional Riemann--Liouville systems 
with  fractional order in $]\frac{2}{3}, 1]$,
under a more general assumption in the nonlinear part, that allows one 
to cover new problems with relevance in some real problems, for example
processes with impulsive type initial condition to describe certain 
characteristics of viscoelastic materials \cite{Du}. 

The text is organized in the following way. In Section~\ref{sec:2}, 
we present the considered problem, we recall some definitions of fractional operators 
involving Riemann--Liouville and Caputo derivatives, and we also recall the definition 
of mild-solution, respectively for Riemann--Liouville fractional systems 
(with a left sided derivative) and Caputo systems with right-sided fractional derivatives. 
Some useful properties are also introduced, as the regional controllability concept 
for the kind of systems under investigation. In Section~\ref{sec:3}, 
after extending the notion of regional controllability to the boundary case,
we establish a relation between regional and boundary regional controllability. 
We prove our main result under an assumption in the nonlinear term that often 
is satisfied in applications of Biology. The result is proved by extending 
the steps of the HUM method to the fractional-order 
case under study and by using Schauder's fixed point theorem. 
To illustrate the obtained result numerically, 
in Section~\ref{sec:4} we apply the proposed method 
to study the controllability on a part of the boundary 
of a Riemann--Liouville fractional version 
of the diffusion logistic growth law model. 
Our example shows the importance of choosing 
an appropriate region $\omega$ and the location 
of the actuator to obtain more efficient results 
on $\varGamma$ with $\varGamma\subset\partial\omega$. 
We end with Section~\ref{sec:5} of conclusions, 
pointing also some possible directions of future research.


\section{Problem Setting and Preliminaries}
\label{sec:2}

We start by recalling some necessary definitions
from fractional calculus and semigroup theory.

\begin{definition}[See \cite{11}]
The left sided Riemann--Liouville integral of order $\alpha>0$ 
for a given integrable function $y$ is defined as
\begin{equation*}
\begin{array}{lll}
\mathrm{I}_{0^+}^\alpha y(t)
&=& \dfrac{1}{\Gamma(\alpha)}\displaystyle\int_{0}^{t}(t-s)^{\alpha-1}y(s)ds, 
\qquad t>0,
\end{array}
\end{equation*}
where $\Gamma$ is the gamma function.
\end{definition}

\begin{definition}[See \cite{11}]
The left sided Riemann--Liouville derivative of order $0<\alpha < 1$ 
is defined as
\begin{equation*}
\begin{array}{lll}
^{RL}\mathrm{D}_{0^+}^{\alpha}
&=&	\dfrac{d}{dt}{\mathrm{I}_{0^+}^{1-\alpha}}.	
\end{array}
\end{equation*}
\end{definition}

\begin{definition}[See \cite{11}]
The right sided Caputo derivative of order $0<\alpha < 1$
for a given differentiable function $y$ is defined as
\begin{equation*}
\begin{array}{lll}
\displaystyle ^C\mathrm{D}_{T^-}^\alpha y(t)
&=& \dfrac{-1}{\Gamma(1-\alpha)}\displaystyle\int_{t}^{T}
(s-t)^{\alpha-1}\dfrac{dy(s)}{ds} ds, 
\qquad  \  0\leq t<T.
\end{array}
\end{equation*}
\end{definition}

\begin{definition}[See \cite{pazy2012semigroups}]
Let $\mathrm{X}$ be a Banach space. We define a $C_0$-semigroup 
(a strongly continuous semigroup) to be a family 
$(S(t))_{t\in\mathbb{R}^+}$ of bounded linear operators 
from $\mathrm{X}$ into $\mathrm{X}$ that satisfies:
\begin{itemize}
\item [i.] $S(0)=I$, where $I$ is the identity operator of $\mathrm{X}$;
\item [ii.] $S(t+s)=S(t)+S(s)$ for all $t,s\in \mathbb{R}^+$;
\item[iii.] $\lim\limits_{t\longrightarrow 0^+} ||S(t)x-x||_X=0$ 
for all $x\in\mathrm{X}$.
\end{itemize}
\end{definition}

\begin{definition}[See \cite{pazy2012semigroups}]
The infinitesimal generator of a $C_0$-semigroup 
$(S(t))_{t\in\mathbb{R}^+}$ is the following linear unbounded operator: 
$$
\mathrm{A}x=\lim\limits_{t\longrightarrow 0^+}\dfrac{S(t)x-x}{t},
$$
defined for every $x$ in its domain given by 
$$
D(\mathrm{A})=\left\{x\in\mathrm{X}  :  
\lim\limits_{t\longrightarrow 0^+}\dfrac{S(t)x-x}{t} 
\ \ \text{ exists } \right\}.
$$ 
\end{definition}

Now we introduce the considered problem in this work.
Let $\Omega$ be an open bounded and regular subset of $\mathbb{R}^n$, $n\leq 3$, 
with a smooth boundary $\partial\Omega$. For $T>0$, we denote $Q=\Omega\times [0,T]$ and 
$\Sigma=\partial\Omega\times [0,T]$. We consider the following semilinear time-fractional
sub-diffusion system with Riemann--Liouville fractional derivative 
of order $\alpha \in \, ]1/2,1]$: 
\begin{equation}
\label{rlcb}	
\left\{
\begin{array}{lllll}
^{RL}\mathrm{D}^\alpha_{0^+}y(x,t)
&=&\mathrm{A}y(x,t)+ Ny(x,t)+\mathcal{B}u(t) & \mbox{in} \quad & Q,\\
\quad \lim\limits_{t\rightarrow 0^+}\mathrm{I}^{1-\alpha}_{0^+}y(x,t)
&=&y_0 & \mbox{in} \quad& \Omega,\\
\dfrac{\partial y(\xi,t)}{\partial \nu_A}
&=&0& \mbox{on} \quad & \Sigma,
\end{array} 
\right. 
\end{equation}
where $^{RL}\mathrm{D}^\alpha_{0^+}$ and $\mathrm{I}^{1-\alpha}_{0^+}$  
denote, respectively, the Riemann--Liouville fractional derivative 
of order $\alpha$ and the Riemann--Liouville fractional integral  of order 
$(1-\alpha)$. The operator $\mathrm{A}$ is an infinitesimal generator of a 
$C_0$-semigroup $(\mathcal{S}(t)_{t>0})$  on the Hilbert space 
$\mathrm{X} = H^1(\Omega)$. In addition, $y\in L^2(0,T;X)$, 
the initial datum $y_0$ is in $\mathrm{X}$, 
$\nu_A$ is the outward unit normal to $\partial\Omega$ 
of the operator  $\mathrm{A}$, and $\mathcal{B}$ is the control
operator, which is a bounded (possibly unbounded) linear operator 
from $\mathbb{R}^{m}$ into $\mathrm{X}$, where $m$ is the number 
of actuators, $u\in \mathcal{U}:=\mathrm{L}^2(0,T; \mathbb{R}^m)$,  
and $N$ is a Lipschitz continuous nonlinear operator.
The Riemann--Liouville fractional semilinear system \eqref{rlcb} 
is defined for $\alpha \in \, ]1/2, 1]$ because this is a condition 
for existence and uniqueness of a mild solution.

We proceed by recalling the necessary results 
and notions to be applied throughout the paper.

\begin{Lemma}[See \cite{11}]
For all $\beta\geq 0$, the  probability density 
function $\mathrm{\phi}_\alpha$ satisfies
$$ 
\int_{0}^{\infty}\theta^\beta \mathrm{\phi}_\alpha (\theta) d\theta
=\dfrac{\Gamma(1+\beta)}{\Gamma(1+\alpha\beta)}.
$$
\end{Lemma}

\begin{Lemma}[See \cite{exis}]
Suppose that $f\in L^2(0,T;\mathrm{X})$ where $\mathrm{X}$ 
is a Banach space. Then, for any initial datum $y_0$, the problem 
\begin{equation}	
\left\{\begin{array}{lllll}
^{RL}\mathrm{D}^\alpha_{0^+}y(t)
&=&\mathrm{A}y(t)+ f(t),\\
\quad \lim\limits_{t\rightarrow 0^+}\mathrm{I}^{1-\alpha}_{0^+}y(t)
&=&y_0 \in D(A), 
\end{array} 
\right. 
\end{equation}
has a unique mild-solution $y\in L^2(0,T;\mathrm{X})$ given by
$$
y(t)=\mathbb{P}_\alpha(t)y_0
+\displaystyle\int_{0}^{t} \mathbb{P}_\alpha(t-s)f(s)ds,
$$
where 
$$
\mathbb{P}_\alpha(t)=\alpha\displaystyle\int_{0}^{\infty}
\theta t^{\alpha-1}\mathrm{\phi}_\alpha(\theta)
\mathcal{S}(t^\alpha\theta)d\theta 
$$
with $\mathrm{\phi}_\alpha$ the probability density function. 
\end{Lemma}

Now, consider the following Caputo fractional differential equation:
\begin{equation}
\label{cap}
\left\{
\begin{array}{lllll}
{^C}D_{T^-}^\alpha p(t) &=& Ap(t),\\
p(T)&=&p_T.
\end{array} \right.
\end{equation}

We consider the following definition of a mild solution for \eqref{cap}.

\begin{definition}[See \cite{mee}]
The function $p(t)\in L^2(0,T;X)$ defined by
$$
p(t)=\mathbb{S}_\alpha(T-t)p_T
$$
is a mild solution of system \eqref{cap},
where 
$$ 
\mathbb{S}_\alpha(t) 
= \displaystyle\int_{0}^{\infty}\mathrm{\phi}_\alpha(\theta)
\mathcal{S}(t^\alpha\theta)d\theta.
$$
\end{definition}

At this point, we can define the mild solution $y(x,t,u):=y_u(t)$ 
of system \eqref{rlcb} by the following integral equation:
\begin{equation}
\label{sys1.sol}
\begin{array}{ll}
y_u(t) = \mathbb{P}_\alpha(t)y_0 
+ \displaystyle\int_{0}^{t}\mathbb{P}_\alpha(t-\tau)
\left(Ny_u(\tau)+\mathcal{B}u(\tau)\right)d\tau.
\end{array} 
\end{equation} 

\begin{Lemma}[See \cite{majo}]
\label{iks}
For $\alpha\in (0,1)$ and $t\in[0,T]$ one has
\begin{equation}
\label{ikss}
\mathrm{I}_{0^+}^{1-\alpha}(\mathbb{P}^*_\alpha(t))
= \mathrm{S}^*_\alpha(t), 
\end{equation}
where $\mathbb{P}^*_\alpha$ and $\mathrm{S}^*_\alpha$ 
are the adjoint operators of  $\mathbb{P}_\alpha$ 
and $\mathrm{S}_\alpha$, respectively.
\end{Lemma}

\begin{Lemma}[See \cite{ren}]
\label{major}
For any $t>0$, $\mathrm{K}_\alpha(t):=t^{1-\alpha}\mathbb{P}_\alpha(t)$ 
is a continuous linear bounded operator, i.e., there exists $M> 0$ such that
$$
\|\mathrm{K}_\alpha(t)\|_{\mathcal{L}(\mathrm{X}, \mathrm{X})}
\leq \dfrac{M}{\Gamma(\alpha)}.
$$
\end{Lemma}

\begin{Lemma}[See \cite{majori}]
\label{inegalit}
For $\gamma\in(0,1] $ and $0<t_1\leq t_2$, the following inequality holds:
$$
|t_1^\gamma-t_2^\gamma|\leq (t_2-t_1)^\gamma.
$$
\end{Lemma}

Let us consider $\omega\subseteq \Omega$ and denote 
the projection operator on $\omega$ by the following map:
$$
\begin{array}{llll}
\chi_{_{\omega}} :& H^1(\Omega) &\longrightarrow & H^1(\omega)\\
&y & \longmapsto &  y_{|_{\omega}}.
\end{array}
$$
We denote by $\chi^*_{_{\omega}}$ its adjoint operator.

We finish this section with the notion
of internal regional controllability.

\begin{definition}[See \cite{ren}]
(i) System \eqref{rlcb} is said to be $\omega$-exactly controllable if for all 
$y_d\in \mathrm{X}$ there exists a control $u\in \mathcal{U}$ such that 
$\chi_{_{\omega}}y_u(T)=y_d$.
(ii) System \eqref{rlcb} is said to be $\omega$-approximately controllable 
if for all $y_d\in ~\mathrm{X}$, given $\varepsilon\geq 0$, there exist a control 
$u\in\mathcal{U}$ such that 
$$
\|\chi_{_{\omega}}y_u(T)-y_d\|_{H^1(\omega)} \leq \varepsilon.
$$
\end{definition}

The main result of our work is to study the regional boundary 
controllability in the semilinear case, which is the subject of 
Section~\ref{sec:3}. To obtain our results for the semilinear case, 
we assume that the linear part of our semilinear system is regionally 
controllable and we use available results in the literature about 
the linear controllability. According to the choice 
of the supports of actuators, for the controllability 
of the associated linear part in a subregion $\omega$  
we can use Proposition~\ref{Prop:2.13}.

\begin{proposition}[See \cite{regana}]
\label{Prop:2.13}
If we denote 
\begin{equation*}
\mathrm{H}^\alpha_\omega u=\displaystyle\int_{0}^{T}
\chi_{\omega}\mathbb{P}_\alpha(T-s)\mathcal{B}u(s)ds,
\end{equation*}
then the following two properties are equivalent:
\begin{enumerate}
\item the linear system associated to \eqref{rlcb} 
is approximately $\omega$-controllable at time $T$;
\item $\overline{Im \mathrm{H}^\alpha_\omega}=H^1(\omega)$.
\end{enumerate}
\end{proposition}


\section{The Boundary Regional Controllability Problem}
\label{sec:3}

In this section, we explore the Hilbert Uniqueness Method (HUM), 
first introduced by Lions in \cite{lions,lion}, and further developed 
in 2016 and 2017 to cover the case of fractional linear systems \cite{art2,art1}. 
Our goal is to find a control that ensures the boundary regional controllability. 
For that, we make use of a relation between internal and boundary regional controllability.

Let us consider the trace operator of order zero $\gamma_0: H^1(\Omega) 
\to H^{\frac{1}{2}}(\partial\Omega)$ such that $\gamma_0 y=y|_{\partial\Omega}$.  
The onto operator $\gamma_0$ is linear and continuous.
Let $\gamma_0^*$ denote its adjoint operator,
which is the solution operator of the Laplace equation
\begin{equation}
\label{ps} 
\left\{
\begin{array}{lllll}
\Delta u
&=&0 & \mbox{in} & \Omega,\\
u+\gamma_0 u&=&0 &\mbox{on}& \partial\Omega.   
\end{array}
\right. 
\end{equation}
The solution of \eqref{ps} exists in general. However, 
an explicit expression for $\gamma_0^*$ is not available.
For more  details and properties about the trace 
and	adjoint operators, we refer the reader to \cite{Touhami}. 

If $\mathrm{\varGamma}\subseteq \partial\Omega$, then we consider 
$\chi_{_{\mathrm{\varGamma}}}$ given by
$$
\begin{array}{rlll}
\chi_{_{\mathrm{\varGamma}}} :
& H^{\frac{1}{2}}(\partial\Omega) 
&\longrightarrow & H^{\frac{1}{2}}(\mathrm{\varGamma})\\
&y & \longmapsto &  y_{|_{\mathrm{\varGamma}}},		
\end{array}
$$
and we denote by $\chi^*_{_{\mathrm{\varGamma}}}$ 
its associated adjoint operator.

\begin{definition}
We say that system \eqref{rlcb} is exactly (respectively approximately) 
boundary controllable on $\mathrm{\Gamma}$ ($\mathcal{B}$-controllable 
on $\mathrm{\Gamma}$) if for all $y_d\in H^{\frac{1}{2}}(\mathrm{\varGamma})$ 
and for all $\varepsilon>0$, there exists a control $u\in\mathcal{U}$ such that
$$
\chi_{_{\mathrm{\varGamma}}}(\gamma_0 y_u(T))=y_d 
\quad  \left(\mbox{respectively } 
\|\chi_{_{\mathrm{\varGamma}}}(\gamma_0
y_u(T))-y_d\|_{H^{\frac{1}{2}}(\mathrm{\varGamma})}
\leq \varepsilon \right). 
$$
\end{definition}

The following lemma gives a useful relation between internal 
and boundary  regional controllability on $\varGamma$.

\begin{Lemma}[See \cite{regana,mee}]
\label{intbound}
If $\mathrm{\varGamma}\subseteq \partial\omega$  
and system \eqref{rlcb} is $\omega$-exactly (respectively $\omega$-approximately) 
controllable, then it is exactly (respectively approximately) $\mathcal{B}$-controllable 
on $\mathrm{\varGamma}$.
\end{Lemma} 

Throughout the rest of the paper we assume that 
$\mathrm{\varGamma}\subseteq \partial\omega$, 
so that it is sufficient, by Lemma~\ref{intbound}, 
to prove the controllability results on a constructive internal part  
$\omega$ of the evolution domain $\Omega$  such that 
$\mathrm{\varGamma}\subseteq \partial\omega$. 

To obtain regional controllability in $\omega$  
of system \eqref{rlcb}, we consider the following hypotheses:
\begin{itemize}
\item [$({H_0})$]  The nonlinear operator $N$ satisfies: 
\begin{equation}
\label{csn}
\exists (L,K)\neq (0,0)\  \mbox{ a positive pair,}  
\  \mbox{ such that } \  
\|Ny\|_X\leq L\|y\|_X+K\|y\|^2_X.
\end{equation}

\item[$({H_1})$] The fractional order derivative satisfies 
$\dfrac{2}{3}< \alpha\leq 1$.
\end{itemize}

The condition $\alpha\in (2/3,1]$ of hypothesis $({H_1})$ 
is necessary in order to deal with the class 
of semi-linear systems that we are studying
(see the proof of Theorem~\ref{thm:01}).
Moreover, the interval $(2/3,1]$ represents a region 
that includes the classical model ($\alpha = 1$) 
and is large enough to describe well the deviation 
of the classical behavior in several applications \cite{bal}.

While system \eqref{rlcb} is well defined for 
$1/2 < \alpha\leq 1$, we are only able to prove	
our regional controllability result under hypotheses $({H_0})$
and $({H_1})$. We can see that $(H_0)$ is a kind of generalization 
of the conditions used in the HUM approach, and that there are 
many nonlinear operators satisfying the $(H_0)$ condition, 
which means that several real phenomena are included in our model.
For example, $({H_0})$ holds for some important
nonlinear operators that appear in biological systems, 
including the logistic growth law model \cite{exp,pop}.

\begin{Theorem}
\label{thm:01}
Assume that $(H_0)$ and $(H_1)$ hold. 
In addition, let the associate linear system to \eqref{rlcb} be approximately 
$\omega$-controllable at time $T$. Then the semilinear system \eqref{rlcb} 
is exactly controllable in $\omega$ by the control function 
$$
u^*(t)=\mathcal{B} ^{*^{RL}}\mathrm{D}_{T^-}^{1-\alpha}\varphi(t),
$$
where $\varphi$ is the mild solution of the following retrograde system:
\begin{equation}
\label{adj}
\left\{
\begin{array}{lllll}
{^C}D_{T^-}^\alpha \varphi(t) 
&=& A^*\varphi(t), \\
\varphi(T)&=&\chi_{_{\omega}}^*\varphi_T:=\varphi_0\in D(A^*).
\end{array} \right.
\end{equation}
\end{Theorem}

In order to prove Theorem~\ref{thm:01}, we recall the steps of the
fractional HUM approach, that transfers our controllability problem 
to a solvability problem (a fixed point problem to obtain  $\varphi_0$).

Let us consider
$$
G=\left\{ y\in \mathrm{X} \quad \mbox{such that}  
\ y=0 \  \mbox{in} \  \Omega\backslash \omega\right\},
$$
in which we define the norm \cite{mee}: 
$$
\|\varphi_0\|_G
=\| \mathcal{B}^ {*^{RL}}\mathrm{D}_{T^-}^{1-\alpha}
\varphi(\cdot) \|_{U}.
$$
Consider system \eqref{rlcb} controlled by 
$u^*(t)= \mathcal{B}^{* \ ^{RL}}\mathrm{D}_{T^-}^{1-\alpha} \varphi(t)$:
\begin{equation}
\label{rlcont}	
\left\{\begin{array}{lllll}
^{RL}\mathrm{D}^\alpha_{0^+}\phi(x,t)
&=&\mathrm{A}\phi(x,t)+ N\phi(x,t)+\mathcal{B}\mathcal{B}^{* \ ^{RL}}
\mathrm{D}_{T^-}^{1-\alpha} \varphi(t) & \mbox{in} \quad & Q,\\
\quad \lim\limits_{t\rightarrow 0^+}\mathrm{I}^{1-\alpha}_{0^+}\phi(x,t)
&=&y_0 & \mbox{in} \quad& \Omega,\\
\dfrac{\partial \phi(\xi,t)}{\partial \nu_A}
&=&0& \mbox{on} \quad & \Sigma,
\end{array} 
\right. 
\end{equation}
where $\varphi(t)=\mathbb{S}_\alpha^*(T-t)\varphi_0$  
is the mild solution of system \eqref{adj}.
We decompose system \eqref{rlcont} into three parts: 
the first being linear, 
\begin{equation}
\label{psrcrl} 
\left\{
\begin{array}{lll}
^{RL}\mathrm{D}_{0^+}^\alpha \phi_0(t)
&=&\mathrm{A}\phi_0(t),\\
\lim\limits_{t\rightarrow 0^+}\mathrm{I}^{1-\alpha}_{0^+}\phi_0(t)
&=&y_0,   
\end{array}
\right. 
\end{equation}
the second linear controlled 
by $u^*$, 
\begin{equation} 
\label{pscontro} 
\qquad \qquad \ \ \   
\left\{
\begin{array}{lll}
^{RL}\mathrm{D}_{0^+}^\alpha \phi_1(t)
&=&\mathrm{A}\phi_1(t)+\mathcal{B}\mathcal{B}^{*^{RL}}
\mathrm{D}_{T^-}^{1-\alpha}\varphi(t),\\
\lim\limits_{t\rightarrow 0^+}\mathrm{I}^{1-\alpha}_{0^+}\phi_1(t)
&=& 0,   
\end{array}
\right. 
\end{equation}
and the third being semilinear,
\begin{equation}
\label{psnonlin}  
\qquad \ \  \qquad  
\left\{
\begin{array}{lll}
^{RL}\mathrm{D}_{0^+}^\alpha \phi_2(t)
&=&\mathrm{A}\phi_2(t)+N(\phi_0 (t)+\phi_1(t)+\phi_2(t)), \\
\lim\limits_{t\rightarrow 0^+}\mathrm{I}^{1-\alpha}_{0^+}\phi_2(t)
&=&0,   
\end{array}
\right. 
\end{equation}
such that $\phi=\phi_0+\phi_1+\phi_2$. 
Let us now define the operators
$$ 
L_o=\chi_\omega^*\chi_\omega(\phi_0(T)), 
\qquad     
C_l \varphi_0=\chi_\omega^*\chi_\omega(\phi_1(T)), 
\qquad \mbox{and} \qquad 
O_N\varphi_0= \chi_\omega^*\chi_\omega(\phi_2(T)).
$$
Thus, the controllability problem reduces to solve 
\begin{equation}
\label{proeq}
C_l \varphi_0= \chi_\omega^* y_d- L_o- O_N \varphi_0.
\end{equation}
Since the associated linear system of \eqref{rlcb} 
is $\omega$-approximately 
controllable, then, by Lemma~3.5 in \cite{regana}, 
$C_l$ is an invertible operator. If we denote 	
\begin{equation}
\label{proeqeq}
\mathcal{K}\varphi_0:= C_l^{-1}\chi_\omega^* y_d
- C_l^{-1}L_o-C_l^{-1} O_N \varphi_0,
\end{equation}
then the equivalent solvability problem is to find 
a fixed point of the operator $\mathcal{K}$.
We are now ready to prove Theorem~\ref{thm:01}.

\begin{proof}
In order to prove Theorem~\ref{thm:01}, we use Schauder's fixed point theorem. 
For that we need to prove that $\mathcal{K}$ defined by \eqref{proeqeq}
maps $B(0,r):=B_r\subseteq G$ into itself for $r>0$ and, moreover, $\mathcal{K}$ 
is a compact operator. According to \eqref{proeqeq}, $\mathcal{K}$ 
is compact if, and only if, $O_N$ is compact.
Let us now show that $O_N(B_r)=\{O_N \varphi_0 
\quad |\quad  \varphi_0\in B_r\}$ is relatively compact.
The fact that 
$$
O_N(B_r)\subset \mathrm{E}_N
:=\{ \chi_\omega^*\chi_\omega(\phi_2(t))\quad  | \varphi_0\in B_r, \ \  t\in[0,T]\},
$$
tell us that it is sufficient to show that $\mathrm{E}_N$ is relatively compact.
Using the Arzel\`a--Ascoli theorem, we can prove the relative compactness 
of $\mathrm{E}_N$ in two steps:
\begin{itemize}
\item [(i)] $\mathrm{E}_N$ is uniformly bounded.
For simplify, we denote $P_\omega=\chi_\omega^*\chi_\omega$ 
and, by the linearity and continuity of $P_\omega$, we have that
$$ 
\exists c_\omega>0 \ \ \mbox{such that} \ 
\| P_\omega(\phi_2(t))\|_{G}
\leq c_\omega \|\phi_2(t)\|_{\mathrm{X}}. 
$$
Moreover, $\phi_2(t)$ is the mild solution of \eqref{psnonlin}, 
which is written as
$$
\phi_2(t)=\displaystyle\int_{0}^{t}
\mathbb{P}_\alpha(t-\tau)N(\phi_0 (\tau)
+\phi_1(\tau)+\phi_2(\tau))d\tau.
$$
From Lemma~\ref{major} and hypothesis ($H_0$), we have 
$$
\begin{array}{lll}
\|\phi_2(t)\|_{\mathrm{X}}
&\leq& \displaystyle\int_0^t\|(t-s)^{\alpha-1}K_\alpha(t-s)N[\phi_0(s)
+\phi_1(s)+\phi_2(s)]\|_{\mathrm{X}}ds\\
& \leq & \dfrac{M}{\Gamma(\alpha)}\displaystyle\int_0^t(t-s)^{\alpha-1}\|
N[\phi_0(s)+\phi_1(s)+\phi_2(s)]\|_{\mathrm{X}}ds\\
& \leq & \dfrac{M L}{\Gamma(\alpha)}
\displaystyle\int_0^t(t-s)^{\alpha-1}\|\phi_0(s)+\phi_1(s)
+\phi_2(s)\|_{\mathrm{X}}ds\\
& & +\dfrac{M K}{\Gamma(\alpha)}\displaystyle\int_0^t(t-s)^{\alpha-1}
\|\phi_0(s)+\phi_1(s)+\phi_2(s)\|^2_{\mathrm{X}}ds,
\end{array}
$$
and, by the same argument used in \cite{mee}, we get that
$$ 
\begin{array}{llll}
\|\phi_2(t)\|_{\mathrm{X}}
&\leq & \dfrac{M L}{\Gamma(\alpha)}\displaystyle\int_0^t(t-s)^{\alpha-1}
\left[\|\phi_0(s)\|_{\mathrm{X}}+\|\phi_1(s)
\|_{\mathrm{X}}+\|\phi_2(s)\|_{\mathrm{X}}\right]ds\\
& & +\dfrac{3M K}{\Gamma(\alpha)}\displaystyle\int_0^t(t-s)^{\alpha-1}
\left[\|\phi_0(s)\|^2_{\mathrm{X}}
+\|\phi_1(s)\|^2_{\mathrm{X}}
+\|\phi_2(s)\|^2_{\mathrm{X}}\right]ds \\
& \leq &  \mathcal{C}(T,\varphi_0)+\dfrac{M L}{\Gamma(\alpha)}
\displaystyle\int_0^t(t-s)^{\alpha-1} \|\phi_2(s)\|_{\mathrm{X}}ds\\ 
& & +\dfrac{3MK}{\Gamma(\alpha)}\displaystyle\int_0^t
(t-s)^{\alpha-1} \|\phi_2(s)\|^2_{\mathrm{X}}ds,
\end{array} 
$$
where 
$$
\begin{array}{lll}
\mathcal{C}(T,\varphi_0)
&= &\dfrac{M^2 L T^{2\alpha-1}}{(\Gamma(\alpha))^2} 
\left[\| y_0\|_{\mathrm{ X}} 
\mathbf{B}(\alpha, \alpha)+\dfrac{T}{\alpha} 
\| \mathcal{B}\|_U\| \varphi_0
\|_{G} \mathbf{B}(\alpha, \alpha+1) \right]\\
& &+ \dfrac{3M^3 K T^{3\alpha-2}}{(\Gamma(\alpha))^3} 
\left[\| y_0\|^2_{\mathrm{ X}}\mathbf{B}(\alpha, 
2\alpha-1)+\dfrac{T^2}{\alpha^2} \| \mathcal{B}\|^2_U\| 
\varphi_0\|^2_{G} \mathbf{B}(\alpha, 2\alpha+1) \right].
\end{array}
$$
Recalling that 
$$
\phi_0(t)=t^{\alpha-1}\mathrm{K}_\alpha(t)y_0,
$$
and 	
$$
\phi_1(t)=\displaystyle\int_{0}^{t}(t-\tau)^{\alpha-1}
\mathrm{K}_\alpha(t-\tau)\mathcal{B}\mathcal{B}^{*^{RL}}
\mathrm{D}_{T^-}^{1-\alpha}\varphi(\tau)d\tau,
$$
we obtain the following inequalities:
$$ 
\begin{array}{llll}
\displaystyle\int_0^t(t-s)^{\alpha-1} \|
\phi_0(s)\|_{\mathrm{X}} ds
&\leq&\dfrac{M}{\Gamma(\alpha)} \| y_0\|_{\mathrm{ X}} 
T^{2\alpha-1} \mathbf{B}(\alpha, \alpha),\\
\displaystyle\int_0^t(t-s)^{\alpha-1} 
\|\phi_0(s)\|^2_{\mathrm{X}} ds
&\leq&\left(\dfrac{M}{\Gamma(\alpha)}\right)^2 
\| y_0\|^2_{\mathrm{ X}} T^{3\alpha-2} 
\mathbf{B}(\alpha, 2\alpha-1),\\
\displaystyle\int_0^t(t-s)^{\alpha-1} 
\|\phi_1(s)\|_{\mathrm{X}} ds
&\leq&\dfrac{M}{\Gamma(1+\alpha)} \|
\mathcal{B}\|_{U}\| \varphi_0\|_{G} 
T^{2\alpha} \mathbf{B}(\alpha, \alpha+1),\\
\displaystyle\int_0^t(t-s)^{\alpha-1} \|\phi_1(s)
\|^2_{\mathrm{X}} ds&\leq&\left(\dfrac{M}{\Gamma(1+\alpha)}\right)^2 
\| \mathcal{B}\|^2_{U}\| \varphi_0\|^2_{G} 
T^{3\alpha} \mathbf{B}(\alpha, 2\alpha)+1,
\end{array}
$$
which are only valid under assumption $({H_1})$.

By applying the generalized Gronwall's lemma 
in \cite{lemgron} with 
$$
f(s)=\dfrac{ML}{\Gamma(\alpha)}(t-s)^{\alpha-1}, 
\quad
g(s)=\dfrac{MK}{\Gamma(\alpha)}(t-s)^{\alpha-1},
$$ 
and $n=2$, we obtain that
$$
\begin{array}{lll}
\| \phi_2(t)\|_{\mathrm{ X}}
&\leq& \dfrac{\mathcal{C}(T,\varphi_0)
\exp\left({\dfrac{MK T^\alpha}{\Gamma(1+\alpha)}}\right)}{1
-\dfrac{3M T^\alpha K}{\Gamma(1+\alpha)}\exp\left({
\dfrac{ML T^\alpha}{\Gamma(1+\alpha)}}\right)
\mathcal{C}(T,\varphi_0)},
\end{array}
$$
such that $T$ satisfies the following condition: 
$$
1-\dfrac{3M T^\alpha K}{\Gamma(1+\alpha)}
\exp\left({\dfrac{ML T^\alpha}{\Gamma(1+\alpha)}}\right)
\mathcal{C}(T,\varphi_0)>0.
$$
It follows that 
$$
\begin{array}{lll}
\|\mathrm{P}(\psi_2(t))\|_{G}
&\leq& \dfrac{c_\omega\mathcal{C}(T,\varphi_0)
\exp\left({\dfrac{MK T^\alpha}{\Gamma(1+\alpha)}}\right)}{1
-\dfrac{3M T^\alpha K}{\Gamma(1+\alpha)}
\exp\left({\dfrac{ML T^\alpha}{\Gamma(1+\alpha)}}\right)
\mathcal{C}(T,\varphi_0)}. 
\end{array}
$$
Thus, $\mathrm{E}_N$ is uniformly bounded.
\item [(ii)] $\mathrm{E}_N$ is equicontinuous.
For $\varphi_0\in\mathrm{E}_N$, let us consider $t_i,t_s$ 
such that $0<t_i< t_s\leq T$. Then, 
$$ 
\|\phi_2(t_s)-\phi_2(t_i)\|_{\mathrm{X}}\leq \mathrm{J}_1
+\mathrm{J}_2+\mathrm{J}_3,
$$
where
$$
\begin{array}{lll}
\mathrm{J}_1&=& \left\|\displaystyle\int_{0}^{t_i}
[(t_s-\tau)^{\alpha-1}-(t_i-\tau)^{\alpha-1}]
\mathrm{K}_\alpha(t_s-\tau)N(\phi_0(\tau)+\phi_1(\tau)
+\phi_2(\tau))d\tau\right\|_{\mathrm{ X}},\\
\mathrm{J}_2&=& \left\|\displaystyle\int_{0}^{t_i}(t_i-\tau)^{\alpha-1}
[\mathrm{K}_\alpha(t_s-\tau)-\mathrm{K}_\alpha(t_i-\tau)]
N(\phi_0(\tau)+\phi_1(\tau)+\phi_2(\tau)) d\tau \right\|_{\mathrm{ X}},\\
\mathrm{J}_3&=& \left\|\displaystyle\int_{t_i}^{t_s}(t_s-\tau)^{\alpha-1}
\mathrm{K}_\alpha(t_s-\tau)N(\phi_0(\tau)+\phi_1(\tau)
+\phi_2(\tau))d\tau\right\|_{\mathrm{ X}}.
\end{array}
$$
The goal in this step is to show that $J_1$, $J_2$ and $J_3$ tend to zero 
as $t_s-t_i\longrightarrow 0$, independently of the choice of $\varphi_0$ 
in $\mathrm{E}_N$, so that we get that $\mathrm{E}_N$ is equicontinuous.
According to Lemmas~\ref{major} and \ref{inegalit} and assumption \eqref{csn}, 
we obtain that
$$
\begin{array}{lll}\mathrm{J}_1
&\leq &\dfrac{M(t_2-t_1)^{1-\alpha}}{\Gamma(\alpha)}
\displaystyle\int_{0}^{t_i}(t_s-\tau)^{\alpha-1}(t_i-\tau)^{\alpha-1}
L(\|\phi_0(\tau) \|_{\mathrm{X}}
+ \|\phi_1(\tau) \|_{\mathrm{X}} + \|\phi_2(\tau)\|_{\mathrm{X}})d\tau\\ 
&+& \dfrac{M(t_2-t_1)^{1-\alpha}}{\Gamma(\alpha)}\displaystyle\int_{0}^{t_i}
(t_s-\tau)^{\alpha-1}(t_i-\tau)^{\alpha-1}K(\|\phi_0(\tau)\|^2_{\mathrm{X}}
+\|\phi_1(\tau)\|^2_{\mathrm{X}}
+\|\phi_2(\tau)\|^2_{\mathrm{X}})d\tau.
\end{array}
$$
Thus, $\underset{t_s-t_i \longrightarrow 0}{\mathrm{J}_1\longrightarrow 0}$.
On the other hand, for $t_i>0$ and $\nu>0$ small enough, we have
$$ 
\begin{array}{lll}
\mathrm{J}_2
&\leq& \underset{\tau\in[0,t_i-\nu]}{\sup}
\|\mathrm{K}_\alpha(t_s-\tau)-\mathrm{K}_\alpha(t_i-\tau)
\|_{\mathcal{L}(\mathrm{X}, \mathrm{X})}
\displaystyle\int_{0}^{t_i-\nu}(t_i-\tau)^{\alpha-1}\|
N(\phi_0+\phi_1+\phi_2)(\tau)\|_{ \mathrm{ X}}d\tau\\ 
& +& \dfrac{2M }{\Gamma(\alpha)}\displaystyle\int_{t_i-\nu}^{t_i} 
(t_i-\tau)^{\alpha-1}\| N((\phi_0+\phi_1+\phi_2)(\tau))
\|_{\mathrm{X}}d\tau\\
&\leq& \underset{\tau\in[0,t_i-\nu]}{\sup}
\|\mathrm{K}_\alpha(t_s-\tau)-\mathrm{K}_\alpha(t_i-\tau)
\|_{\mathcal{L}(\mathrm{X}, \mathrm{X})}\dfrac{(t_i-\nu)^{\alpha
-\frac{1}{2}}}{(2\alpha-1)^\frac{1}{2}}\\ 
& \times& ( L\| (\phi_0+\phi_1+\phi_2)(\tau)
\|_{ L^2(0,T;\mathrm{X})}
+ K \|(\phi_0+\phi_1+\phi_2)(\tau) \|^2_{ L^4(0,T;\mathrm{X})})\\
&+& \dfrac{2M \nu^{\alpha-\frac{1}{2}}}{(2\alpha-1)^{
\frac{1}{2}}\Gamma(\alpha)}(L \| (\phi_0+\phi_1+\phi_2)(\tau)
\|_{ L^2(0,T;\mathrm{X})}+K \|
(\phi_0+\phi_1+\phi_2)(\tau)) \|^2_{ L^4(0,T;\mathrm{X})}.
\end{array}
$$
Therefore, by the continuity of the operator $\mathrm{K}_\alpha$, 
it can be seen that $\mathrm{J}_2$ tends to zero as 
$t_s-t_i\longrightarrow 0$, $\nu\longrightarrow 0$ 
for all $\varphi_0$ in $\mathrm{E}_N$.
For $\mathrm{J}_3$, we can also obtain the following inequality, 
which ensures that the limit of $\mathrm{J}_3$ 
when $t_s-t_i$ tends to zero is zero:
$$
\mathrm{J}_3
\leq \dfrac{M(t_s-t_i)^{\alpha-\frac{1}{2}}}{(2\alpha-1)^{\frac{1}{2}}
\Gamma(\alpha)}(L \| (\phi_0+\phi_1+\phi_2)(\tau)
\|_{ L^2(0,T;\mathrm{X})}+K \| (\phi_0+\phi_1
+\phi_2)(\tau)\|^2_{ L^4(0,T;\mathrm{X})}). 
$$	 
\end{itemize}
Finally, we prove that $\mathcal{K}(B_r)\subseteq B_r$. 
If this is not the case, then, for each $r>0$, it is possible to find a  
$\varphi_0\in B_r$ such that $\|\mathcal{K}(\varphi_0)\|_{G}> r$.
From the definition of $\mathcal{K}$, we have that
$$
\begin{array}{lll}
r <	\|\mathcal{K}(\varphi_0)\|_{G}
&\leq & \| C_l^{-1}\chi_\omega^* y_d
- C_l^{-1}L_o \|_{G} +  \| C_l^{-1}O_N \varphi_0 \|_{G}\\
& \leq & \| C_l^{-1}\chi_\omega^* y_d- C_l^{-1}L_o \|_{G} 
+ \dfrac{c_\omega\mathcal{C}(T,r) \| C_l^{-1} \|_{{\mathcal{L}(G, G)}}
\exp\left(\dfrac{MK T^\alpha}{\Gamma(1+\alpha)}\right)}{1
-\dfrac{3M T^\alpha K}{\Gamma(1+\alpha)}\exp\left({
\dfrac{ML T^\alpha}{\Gamma(1+\alpha)}}\right)\mathcal{C}(T,r)}. 
\end{array}
$$
Dividing both sides by $r$, and taking the lower limit as $r\rightarrow +\infty$, 
we obtain a contradiction. We conclude that $\mathcal{K}(B_r)\subseteq B_r$.
\end{proof}

It should be noted that Theorem~\ref{thm:01} is not a consequence 
of known results and that its proof employs a novel method. Here we study 
the boundary regional controllability for a class of Riemann--Liouville fractional 
semilinear systems under a condition on the nonlinear part of the system 
that appears in some real models, for example in the nonlinear 
growth population model, that is not covered by the results of
\cite{mee,ma}. To do so, we use in the proof of Theorem~\ref{thm:01} 
the HUM approach and a generalization of the Gronwall--Bellman inequality, 
which contrasts with the method used in papers \cite{mee,ma} 
that is based on the standard Gronwall's lemma.


\section{Algorithm and Simulations}
\label{sec:4}

Now we consider the two-dimensional diffusive logistic population growth nonlinear 
model with a Riemann--Liouville time derivative of fractional order $0<\alpha < 1$:
\begin{equation}
\label{logis}
\left\{
\begin{array}{llll}
^{^{RL}}\mathrm{D}_{0^+}^{^{\alpha}}  N(x,y,t) 
&=& \Delta N(x,y,t) +C N(x,y,t)\left(1-\dfrac{N(x,y,t)}{K}\right) +Bu(t)  
& \text{ in } \Omega\times\left]0,2\right],\\
\dfrac{\partial N(\xi,\mu,t)}{\partial \nu}
&=&0 \qquad  &\mbox{on} \quad  \partial\Omega\times\left]0,2\right],\\
\lim\limits_{t\rightarrow 0^+}\mathcal{I}^{1-\alpha}_{0^+}N(x,y,t) 
&=& N_0(x,y) & \text{ in } \Omega,
\end{array}  \right.
\end{equation} 
where $\Omega=]0,\pi[\times]0,\pi[$, $N$ is the population density, 
$\partial N\slash \partial \nu= \nu\cdot \nabla N|_{\partial\Omega}$,
$t$ is in $[0,2]$, $C$ and $K$ are  positive constants, and where
we consider the initial data $N_0$ to be non negative and nonzero.

Note that the Neumann boundary conditions signify that the number 
of individuals is fixed in the domain. Moreover, the logistic growth 
operator satisfies our hypothesis. In this case, the operators 
$K_\alpha$ and $S_\alpha$ have the following form:
$$ 
\mathrm{S}_\alpha(t)y=\sum_{j=1}^{+\infty} 
\sum_{k=1}^{r_j} \mathrm{E}_{\alpha}(\lambda_j t^\alpha)<y,
\xi_{jk}>_X \xi_{jk}(x) 
\ \mbox{and}  \ \  
\mathrm{K}_\alpha(t)y=\sum_{j=1}^{+\infty} \sum_{k=1}^{r_j} 
\mathrm{E}_{\alpha,\alpha}(\lambda_j t^\alpha)<y,\xi_{jk}>_X \xi_{jk}(x),
$$
where $(\xi_{jk})_{j\geq 1,  1\leq k\leq r_j}$  is a complete system 
of eigenfunctions in $X$ of the Laplace  operator and $\lambda_j$ 
are the associated eigenvalues with the  multiplicity $r_j$.
\begin{algorithm}[H]
\caption{}
\label{algo}
\begin{algorithmic}
\STATE	\textbf{Initialization}:
\STATE  -- Fractional order of derivative $\alpha$, time $T$, 
initial data $y_0$, $\varGamma$ as a part of the boundary $\partial\Omega$, 
the control operator  $B$ (zonal or pointwise), and the desired state $y_d$ on $\varGamma$.
\STATE -- Choose a region $\omega$ such that $\varGamma\subseteq \partial\omega$ 
and an extension of $y_d$ defined in $\omega$. 
\STATE -- Choose  $\varepsilon$ small enough.
\STATE -- Choose an initial function  $\varphi_0^1$.
\STATE -- Let $k = 1$.
\STATE	\textbf{Repeat}
\STATE -- Solve \eqref{adj} and obtain $\varphi^k$. 
\STATE -- Solve \eqref{psrcrl} and obtain $\phi_0^k$.
\STATE -- Solve \eqref{pscontro} and obtain $\phi_1^k$.
\STATE -- Solve \eqref{psnonlin} and obtain $\phi_2^k$.
\STATE -- Put $\varphi_0^{k+1}=\mathcal{K}(\varphi_0^k)$, $k := k+1$.
\STATE	\textbf{Until}
\STATE  \bf{the control} is  
$u^*(t)= \mathcal{B}^{*^{RL}}\mathrm{D}_{T^-}^{1-\alpha}\varphi(t)$.
\end{algorithmic}
\end{algorithm}

By applying Algorithm~\ref{algo}, we show, through a simulation study 
of the diffusion logistic population growth law model, 
the effectiveness of our proposed method to reach the steady state 
at time $T$ by means of a suitable control function 
(zonal or pointwise actuator).

We can calculate the expression of $\varphi$ using the formula
$\varphi(t)=S_\alpha^*(T-t)\varphi_0$. Consequently, after the 
decomposition of the mild solution $N(x,y,t)$, we obtain the expression 
of $\psi_0(t)$ using the formula of $K_\alpha$. For obtaining $\psi_1(t)$, 
we use the expression of $K_\alpha$ and $\mathrm{B}^*$ according 
to the choice of the actuators, while for $\psi_2(t)$ we solve 
the system using the predictor-corrector method and the explicit 
expression of $\psi_0$ and $\psi_1$.


\subsection{Example~1: Zonal Actuator} 
\label{subsect:1}

We consider system \eqref{logis} with $\alpha=0.75$, 
$B=\chi_{_{D}}$,$C=1$, $K=100$ and $N_0(x,y)=\sqrt{xy}$.
The domain of the actuator is $D=[0.5, 1]\times[0.7, 1]$, 
the considered subregion being $\omega=[0, 0.3]\times[0, 0.5]$, 
and the boundary subregion in $\partial\omega$ is 
$\varGamma=\{0\}\times[0, 0.5]$.

\begin{figure}[H]
\begin{subfigure}[b]{0.49 \textwidth}
\begin{center}
\includegraphics[width=6cm]{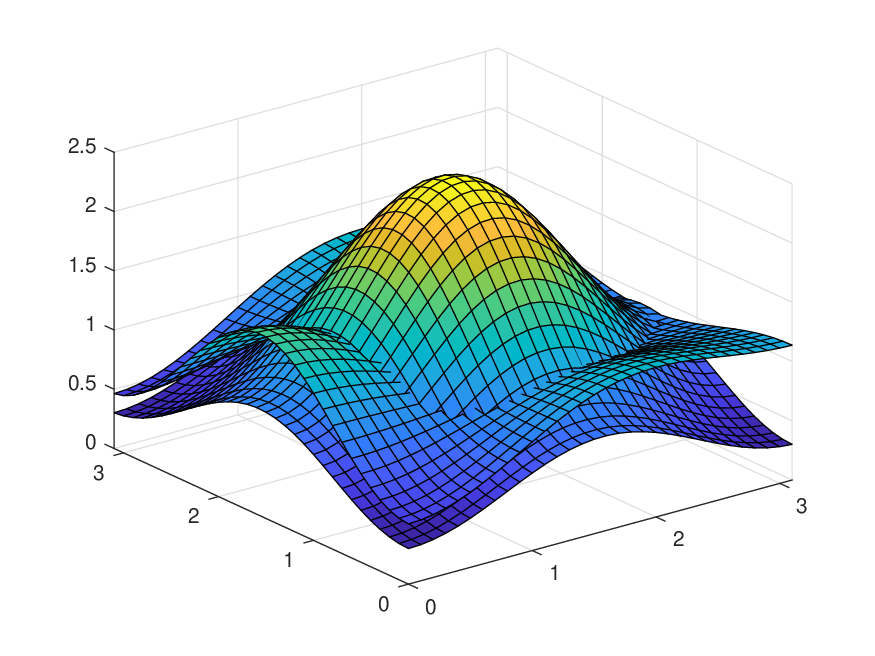}
\caption{The desired and reached states in $\Omega$
for the example of Section~\ref{subsect:1}.}
\label{fig1}
\end{center}
\end{subfigure}
\hfill
\begin{subfigure}[b]{0.5 \textwidth}
\begin{center}
\includegraphics[width=6cm]{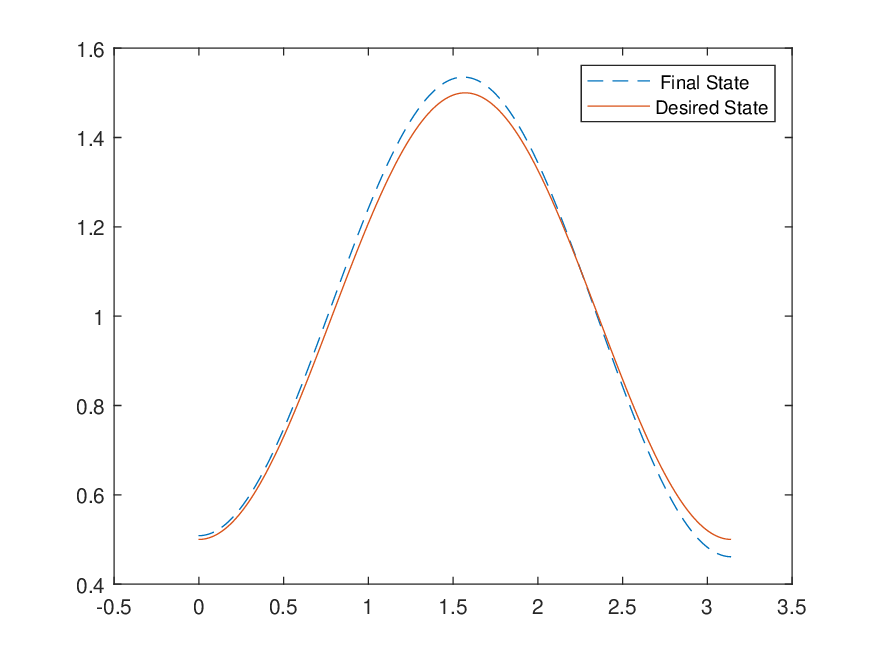}
\caption{The desired and reached states on $\varGamma$
for the example of Section~\ref{subsect:1}.}
\label{fig2}
\end{center}
\end{subfigure}
\end{figure}

Let us consider the desired state $z_d(y)$ on $\varGamma$, 
which is the steady state $H(y)$ of the logistic growth model 
$z_d(y):=H(y)=(1-\mu \cos(2y))$ and let 
$N_d(x,y)=(1-\mu \cos(2x))(1-\delta \cos(2y))$ be the extension 
of the desired state in $\omega$.
By choosing $\mu=0.5$, we obtained Figures~\ref{fig1} and \ref{fig2}.
Precisely, Figures~\ref{fig1} and \ref{fig2} display the desired  
and the reached states in $\omega$ and also on $\varGamma$.
They show that the reached state is very closed to the desired state with an error 
in $\omega$ equal to $8.0\times 10^{-3}$ and an error of order $10^{-5}$ on $\varGamma$.
Our simulation results, including the error on $\varGamma$ as a function of the chosen 
subregion and the actuator location, are summarized in Table~\ref{table:1}, 
which shows that the choice of the region $\omega$ is crucial
to obtain the boundary regional controllability.
\begin{table}[H]
\centering
\caption{Relation Error--Actuator--Region for the 
example of Section~\ref{subsect:1}.}\label{table:1}
\begin{tabular}{|c|c|c|} \hline 
\textbf{Actuator} & \textbf{Region $\mathbf\omega$} 
& \textbf{Error on $\mathbf\varGamma$} \tabularnewline \hline 
$[0.5, 1]\times[0.7, 1]$ & \multirow{2}{1cm}{} 
$[0, 0.5]\times[0, 0.7]$ & $0.1179$ \tabularnewline
\cline{2-3} 
& $[0, 1]\times[0, 1]$ & $0.6392$   \tabularnewline
\hline 
{$[0.3, 0.5]\times[0.7, 1]$} & \multirow{2}{1.1cm}{} 
$[0, 0.3]\times[0, 0.5]$ & $0.0554$ \tabularnewline
\cline{2-3} 
& $[0, 0.1]\times[0, 0.7]$ & $0.0069$   \tabularnewline
\hline
{$[0, 0.3]\times[0, 0.1]$} & 	\multirow{2}{1.1cm}{} 
$[0, 0.5]\times[0, 0.5]$ & $0.0953$ \tabularnewline
\cline{2-3} 
& $[0, 0.5]\times[0, 0.7]$ & $0.1081$   \tabularnewline
\hline
\end{tabular}
\end{table}


\subsection{Example~2: Pointwise Actuator} 
\label{subsect:2}

In this example, we fix $\alpha=0.8$, $B=\delta(b_1,b_2)$ 
(pointwise actuator), $C=1$, $K=1$, and 
$$
N_0=\sqrt{xy}\exp(xy). 
$$
Let us consider $(b_1,b_2)=(0, 0.5)$, 
the desired state 
$$
z_d(y)=100\left(1-\dfrac{1}{2} \cos(2y)\right)
$$ 
on $\varGamma=\{0\}\times[0, 0.4]$ 
and let $N_d(x,y)=2\left(1-\dfrac{1}{2} 
\cos(2x)\right)\left(1-\dfrac{1}{2} \cos(2y)\right)$ 
be the extension of the desired state 
in $\omega=[0, 0.3]\times[0, 0.4]$.
\begin{figure}[H]
\begin{subfigure}[b]{0.49 \textwidth}
\begin{center}
\includegraphics[width=6cm]{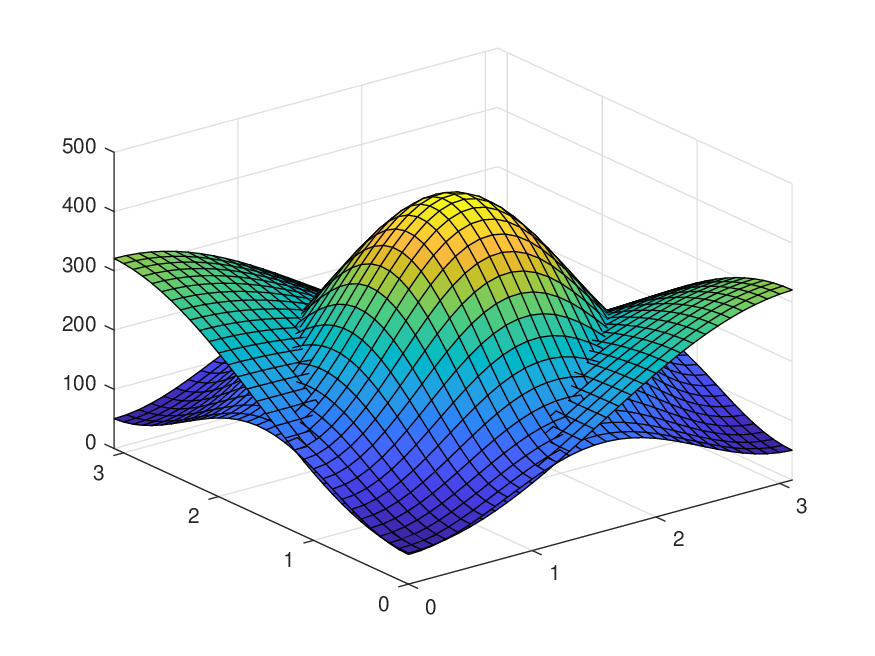}
\caption{The desired and reached states in $\Omega$ 
for the example of Section~\ref{subsect:2}.}
\label{fig12}
\end{center}
\end{subfigure}
\hfill
\begin{subfigure}[b]{0.49 \textwidth}
\begin{center}
\includegraphics[width=6cm]{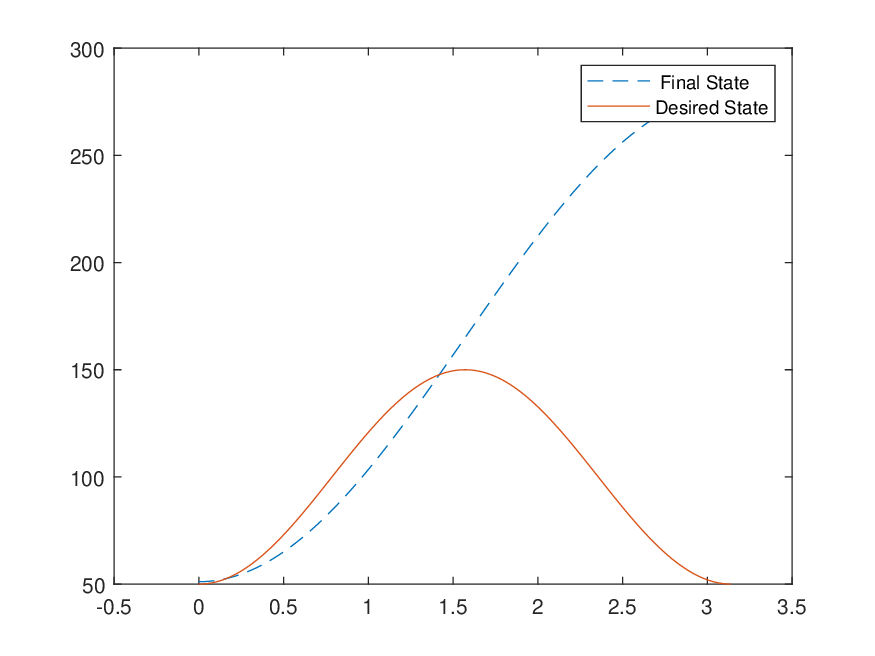}
\caption{The desired and reached states on $\varGamma$ 
for the example of Section~\ref{subsect:2}.}
\label{fig22}
\end{center}
\end{subfigure}
\end{figure}

The results obtained by the application of our method
are reported in Figures~\ref{fig12} and \ref{fig22}.  
They show that the desired state and the steady state 
are close in $\omega$, differing by $9.31\times 10^{-2}$, 
and on $\varGamma$ by only $6.95\times 10^{-5}$.
\begin{table}[H]
\centering
\caption{Relation Error--Actuator--Region for the 
example of Section~\ref{subsect:2}.}\label{table:2}
\begin{tabular}{|c|c|c|}
\hline 
\textbf{Region $\mathbf\omega$} & \textbf{Pointwise actuator} 
& \textbf{Error on $\mathbf \varGamma$} \tabularnewline \hline 
$[0, 0.3]\times[0, 0.4]$ & \multirow{2}{1cm}{} $(0, 1)$ 
& $0.3884$ \tabularnewline \cline{2-3} 
& $(0, 0.3)$ & $0.1022$   \tabularnewline
\hline 
{$[0, 0.5]\times[0, 1]$} 
& \multirow{2}{1.1cm}{} $(0, 0.5)$ 
& $0.3969$ \tabularnewline \cline{2-3} 
& $(0, 0.1)$ & $0.0639$   \tabularnewline \hline
{$[0, 1]\times[0, 0.5]$} & 	\multirow{2}{1.1cm}{} $(0, 0.5)$ 
& $0.0363$ \tabularnewline \cline{2-3} 
& $(0, 0.3)$ & $0.0552$   \tabularnewline \hline
\end{tabular}
\end{table}
We observe in Table~\ref{table:2} that the proposed method 
is highly sensitive with respect to the region $\omega$ 
such that $\varGamma\subseteq \partial\omega$, as well 
as the actuator location.


\section{Conclusion}
\label{sec:5}

The Fractional Hilbert Uniqueness Method (FHUM) has been proposed 
for solving a boundary regional controllability problem for a new
class of semilinear Riemann--Liouville  fractional systems. We studied the 
regional controllability on a subregion of the boundary evolution domain 
under some assumptions on the nonlinear operator and the order of the derivative.  
We developed a numerical simulation to illustrate the obtained result through 
an example of the logistic growth law model with different actuators. An 
interesting open question concerns the study of the gradient regional 
controllability of semilinear fractional systems for fractional operators 
with respect to another function. This is under investigation and will 
be the objective of our future work.


\section*{Acknowledgments}

This research is part of first author's post-doc, which is carried out 
at University of Aveiro, Portugal, within the project ``Mathematical Modelling 
of Multiscale Control Systems: Applications to Human Diseases'' (CoSysM3), 
Reference 2022.03091.PTDC, financially supported by national funds (OE) 
through FCT/MCTES. Tajani and Torres are also supported by Portuguese funds 
through FCT-Portugal within project UIDB/04106/2020 (CIDMA) with DOI  
\url{https://doi.org/10.54499/UIDB/04106/2020}.
The authors are very grateful to four anonymous referees 
and to an Associate Editor for their suggestions and invaluable comments,
which helped them to improve the manuscript.



\end{document}